\documentclass{amsart}

\usepackage{
amsmath,amsthm,amsfonts,amscd,graphicx,latexsym,amssymb,eucal,mathrsfs
}

\usepackage[all]{xy}

\newcommand{\absolute}[1]{\lvert#1\rvert}
\renewcommand{\phi}{\varphi}
\renewcommand{\epsilon}{\varepsilon}
\newcommand{\abelian}{{\rm ab}}

\DeclareMathOperator{\Trace}{Tr}

\DeclareMathOperator{\lcm}{lcm}
\DeclareMathOperator{\dist}{dist}
\DeclareMathOperator{\PGL}{PGL}

\DeclareMathOperator{\branchlocus}{br}

\newcommand{\mathset}[1]{{\{#1\}}}

\newtheorem{thm}{Theorem}[section]
\newtheorem{prop}[thm]{Proposition}
\newtheorem{lem}[thm]{Lemma}
\newtheorem{cor}[thm]{Corollary}
\newtheorem{rem}[thm]{Remark}
\newtheorem{dfn}[thm]{Definition}

\newtheorem{exa}[thm]{Example}

\newtheorem{conv}[thm]{Convention}

\newtheorem{Def}[thm]{Definition}

\newtheorem{Rem}[thm]{Remark}
\newtheorem{Thm}[thm]{Theorem}

\newtheorem{Prop}[thm]{Proposition}

\newtheorem*{Satz*}{Satz}

\newtheorem{Lemma}[thm]{Lemma}

\begin{document}
\title{Cyclic coverings of the $p$-adic projective line by Mumford curves}

\author{Patrick Erik Bradley}

\date{\today}



\begin{abstract}
Exact bounds for the  positions of the branch points for cyclic coverings of the $p$-adic projective line by Mumford
curves are calculated in two ways. Firstly, by using Fumiharu Kato's $*$-trees, and secondly
by giving explicit matrix representations of  the Schottky groups corresponding to the
Mumford curves above the projective line through combinatorial group theory. 
\end{abstract}

\maketitle

\section{Introduction}

Cyclic covers of the projective line defined over a field $K$ of characteristic zero have
been thoroughly studied. Such covers $\phi\colon X\to\mathbb{P}^1$ correspond to equations
of the form 
$$
y^n=f(x),
$$
where $f(x)\in K[x]$ is a polynomial. The zeros of $f(x)$ in a suitable finite extension field 
of $K$ are the branch points of $\phi$.
In the case that $K$ is a $p$-adic field, it is known that  not every equation as above corresponds
to a cover by a Mumford curve. And even if  $f(x)$ is of the right kind, one finds
strong restrictions on the position of the branch points for $\phi$ to be a Mumford cover of $\mathbb{P}^1$.
This was first observed in the case $n=2$ and $X$ an elliptic curve: $X$ is a Tate curve, if and only if the four branch points do not form an equilateral quadrangle in $\mathbb{P}^1$.
To be more precise, by a projective automorphism one can take the branch locus to be $\{0,1,\infty,\lambda\}$ with $\absolute{\lambda}=1$. Then for residue characteristic not equal to two, the Legendre equation
$$
y^2=x(x-1)(x-\lambda)
$$
is the equation of a Tate curve, if and only if $\absolute{\lambda-1}<1$.

If $\phi$ is a hyperelliptic cover, then the restriction found in \cite{Ulrich1981} is that the branch points
come in pairs of points closer to one another than
to the other  branch points. The distance is measured by rational affinoid subsets of $\mathbb{P}^1$. 

\medskip
The most elegant way of obtaining bounds for  relative positions of the branch points of any finite Galois cover $\phi$ of $\mathbb{P}^1$
is through  the $*$-tree $\mathscr{T}^*_N$ for the discrete finitely generated group $N$ giving rise to the orbifold
uniformisation of $\Omega\stackrel{N}{\longrightarrow}\mathbb{P}^1$ factoring through $\phi$ and having the same branch locus and ramification orders
as $\phi$. Distances within $\mathscr{T}_N^*$ translate into distances between branch points.
The group $N$ sits in an exact sequence
$$
1\to \Gamma\to N\to G\to 1,
$$
where $G$ is the Galois group of $\phi$ and $\Gamma$ is a free group whose rank is the genus of $X$.
The corresponding uniformisation $\Omega\stackrel{\Gamma}{\longrightarrow} X$ is called a {\em Schottky uniformisation}, and $\Gamma$
a {\em Schottky group}. This approach is pursued in the present article for cyclic covers.

\medskip
The $*$-tree $\mathscr{T}_N^*$ was developped by Fumiharu Kato in order to obtain deeper insight into the structure of $p$-adic discrete groups
and has  many applications in the study of automorphisms of Mumford curves, especially in positive 
characteristic,
e.g.\ \cite{CKK2001}. 
Extensive use of the $*$-tree is being made in the classification of $p$-adic  triangle
groups \cite{BKV01}. 

\medskip
In the present article, we focus on all possible cyclic covers of $\mathbb{P}^1$ with twofold aim.

\medskip
Firstly, we exhibit detailed calculations of the exact bound for $\absolute{\lambda-1}$ 
characterising the Mumford covers among four-point
cyclic Harbater-Mumford covers $\phi$ whose branch locus is  $\{0,1,\infty,\lambda\}$ and $\absolute{\lambda}=1$,
and give the exact sizes of the separating annuli for covers with more branch points.

\smallskip
Secondly, explicit hyperbolic generators for $\Gamma$ are given from which again one can calculate
the characteristic bound from above by Ford's method of isometric circles and some combinatorial group theory, and thus gains an explicit parametric description of the Schottky uniformisation of the Mumford curve.
For covers of prime degree, we recover the generators of \cite{vanSteen1982}.

We  remark, however,  that here we are only dealing with the positions of branch points
up to ``first order", meaning that our methods do not reveal the precise relationship
between the discrete representation $\Gamma\to\PGL_2(K)$ and the
branch points of the corresponding Mumford cover $\phi$, which would require the study of
automorphic forms on the Mumford curves. Geometrically speaking, we can 
make explicit the geometry of $\mathscr{T}_N^*$ without, however, considering the precise
embeddings of $\mathscr{T}_N^*$ into the Bruhat-Tits tree for $\PGL_2(K)$.

\medskip
A nice desideratum would be  the explicit  fuchsian differential
equation corresponding to the cyclic cover $\phi\colon X\to\mathbb{P}^1$.

\medskip
This article refines methods and results from the author's  dissertation \cite{BradDiss2002}.

\section{Generalities}
Let $\mathbb{Q}_p$ be the field of $p$-adic numbers.
We assume that $K$ is a finite extension field of $\mathbb{Q}_p$, large enough that all branch points 
of all  covers $X\to\mathbb{P}^1_K$ in the article are $K$-rational.

$\mathscr{T}_K$ denotes the Bruhat-Tits tree for $\PGL_2(K)$, the automorphism group of the projective line $\mathbb{P}^1_K$. We will use the well known fact that the ends of $\mathscr{T}_K$ correspond to the $K$-rational points of the projective line $\mathbb{P}^1$. 

Let $N\subseteq\PGL_2(K)$ be a finitely generated discrete subgroup. Following \cite{KatoJAG2005},
the tree $\mathscr{T}_N^*$ is defined to be the 
smallest subtree of $\mathscr{T}_K$ whose ends correspond to the fixed points of all non-trivial elements of $N$.
The group $N$ acts on $\mathscr{T}_N^*$ without inversion, and the quotient graph $T^*_N=\mathscr{T}_N^*/N$ is a graph of finite groups with finitely many ends corresponding to the branch points of the quotient cover
$$
\Omega_N\stackrel{N}{\longrightarrow} X_N.
$$
The open  analytic space $\Omega_N\subseteq \mathbb{P}^1$ is defined as the complement of the closure of the limit points of $N$, and the quotient space $X_N$ is the analytification of a non-singular projective algebraic curve over $K$. In fact, $X_N$ is a Mumford curve.

\smallskip
An important example of $\mathscr{T}^*_N$ is, when $N=\langle\gamma\rangle\cong C_m$ is a finite cyclic group
of order $m>1$. Then $M(\gamma):=\mathscr{T}^*_N$ is simply a straight line stabilised by $\gamma$.

\begin{Def}
Let $\gamma\in\PGL_2(K)$ be of finite positive order. Then $M(\gamma)$ is called the {\em mirror} of $\gamma$.
\end{Def}

\begin{lem}	\label{mirrormaxcyc}
There is a natural bijection between the sets:
$$
\left\{
\text{maximal finite cyclic subgroups of $N$}
\right\}
\stackrel{\sim}{\longrightarrow}
\{\text{mirrors of $N$}\}
$$
\end{lem}

\begin{proof}
The natural map takes a maximal cyclic group to the mirror of a generator, which is clearly well defined.
It is also clear that  $\langle\gamma\rangle\subseteq\langle\delta\rangle$
implies $M(\gamma)=M(\delta)$. Therefore, the map is surjective.

Let now $\langle\gamma\rangle, \langle\delta\rangle\subseteq N$ be such that the mirrors $M:=M(\gamma)=M(\delta)$
coincide. Then $G:=\langle\gamma,\delta\rangle\subseteq N$ is finite,
as any word in $\gamma$ and $\delta$ fixes the mirror. This means that $M=\mathscr{T}^*_G$, implying that
$G$ is cyclic, as the corresponding cover $\mathbb{P}^1\stackrel{G}{\longrightarrow}\mathbb{P}^1$ has
exactly two branch points \cite[Proposition 5.6.2]{KatoJAG2005}. 
\end{proof}

It is well known that, if $K'/K $ is a finite field extension, then a subdivision of $\mathscr{T}_K$
embeds into $\mathscr{T}_{K'}$. 

\begin{conv}
Let $\mathscr{T}$ be a subtree of $\mathscr{T}_K$. When we speak of a {\em point} $x$ on an edge $e=(v,w)$ of $\mathscr{T}$,
we mean that after some finite extension $K'/K$,  $x$ is a vertex on the (open) {\em path} $(v,w)$ in the tree
$\mathscr{T}'$ obtained by restricting the  subdivision and embedding process from above.
\end{conv}

\section{Mumford curves and discrete groups} \label{mumfdiscr}

\subsection{The tree $\mathscr{T}^*_\Gamma$ for a free product of cyclic groups}

\begin{figure}[h]
\begin{center}
\setlength{\unitlength}{1cm}
\begin{picture}(8.5,2.5)
\scriptsize
\thicklines
\put(2.45,1.2){$x$}
\put(3.55,1.2){$v$}
\put(4.3,1.2){$w$}
\put(5.45,1.2){$y$}
\put(2.5,1){\circle*{.15}}
\put(2.5,1){\line(1,0){.3}}
\put(3.6,1){\line(1,0){.8}}
\put(5.2,1){\line(1,0){.3}}
\put(5.5,1){\circle*{.15}}
\put(1.9,.9){$C_{m}$}
\put(5.7,.9){$C_{n}$}
\put(3.95,1.1){$1$}
\thinlines
\put(2.5,1){\vector(-1,1){1}}
\put(2.5,1){\vector(-1,-1){1}}
\put(5.5,1){\vector(1,1){1}}
\put(5.5,1){\vector(1,-1){1}}
\thicklines
\put(4.4,1){\circle*{.15}}
\put(3.6,1){\circle*{.15}}
\put(3.5,.7){$C_p$}
\put(4.2,.7){$C_p$}
\put(2.8,1){\circle*{.15}}
\put(5.2,1){\circle*{.15}}
\put(2.6,.7){$C_{p^{r}}$}
\put(5,.7){$C_{p^{s}}$}
\put(2.8,1){\dashbox{.1}(.8,0){}}
\put(4.4,1){\dashbox{.1}(.8,0){}}
\put(1.2,2.1){$m$}
\put(1.2,-.2){$m$}
\put(6.6,2.1){$n$}
\put(6.6,-.2){$n$}
\end{picture}
\end{center}
\caption{The quotient $*$-tree for a free product of cyclic groups.}
\label{generalsegment}
\end{figure}
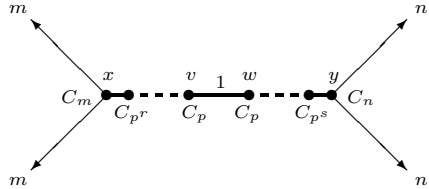

Let $\Gamma=C_m*C_n$ be the free product of two cyclic groups $C_m$ and $C_n$. We will calculate the 
tree $\mathscr{T}_\Gamma^*$ for all possible values of $m=p^ra$ and $n=p^sb$ (where $(a,p)=(b,p)=1$). In fact, the shape of the quotient tree $T_{\Gamma}^*$ is
known in \cite[\S 8.1]{KatoJAG2005} (and implicitly known in \cite[\S 11]{Her78}) and is given in Figure \ref{generalsegment}. 

Our special interest lies in the exact distances within the tree, in particular the lengths of the paths $[x,v]$ and $[w,y]$. These can be extracted from the indications at the end of 
\cite[\S 8.1]{KatoJAG2005} or the proof of \cite[Proposition 3.1]{BradMZ2005}. Here, we give a detailed
exposition of the calculations.

\begin{prop}	\label{*tree4ptcycp'p'}
Let $\Gamma\cong C_m*C_n$ be a discrete subgroup of $\PGL_2(K)$. If $(m,p)=(n,p)=1$, then
$T_\Gamma^*$ is as in Figure \ref{generalsegment} with $\dist(x,v)=\dist(w,y)=0$. 
\end{prop}

\begin{proof}
\cite[\S 8.1]{KatoJAG2005}.
\end{proof}

\begin{prop} \label{*tree4ptcycpp}
Let $\Gamma_1\subseteq\Gamma_2$ be discrete subgroups of $\PGL_2(K)$, where as abstract groups
$\Gamma_1\cong C_p*C_p$ and $\Gamma_2\cong C_{pa}*C_{pb}$ with $a,b\ge 1$. Then: 
\begin{enumerate}
\item There is a subdivision
$\mathscr{T}_1^*$ of $\mathscr{T}_{\Gamma_1}^*$ which is a subtree of $\mathscr{T}_{\Gamma_2}^*$.  
\item The quotient graphs $T_{\Gamma_1}^*$ and $T_{\Gamma_2}^*$ are trees with shape of that in Figure
\ref{generalsegment}.
\item For a primitive $p$-th root $\zeta_p$ of unity, $$\dist(x,v)=\dist(w,y)=v(\zeta_p-1)$$
in both trees $T_{\Gamma_1}^*$ and $T_{\Gamma_2}^*$.  
\end{enumerate}
\end{prop}

\begin{proof}
Since $\Gamma_2$ is a free product of non-trivial cyclic groups, $\mathscr{T}_{\Gamma_2}^*$ contains an edge $e=(v,w)$ with stabilisers $\Gamma_{2,e}=1$ and $\Gamma_{2,v}$, $\Gamma_{2,w}$ both non-trivial. 
The group $\Gamma_{2,v}=\langle\gamma\rangle$ contains an element of order $p$, as otherwise 
$v$ would lie on the mirror $M(\gamma)$ which in turn corresponds to some maximal finite cyclic subgroup of
$\Gamma_2$ (Lemma \ref{mirrormaxcyc}). But such subgroups necessarily contain elements of order $p$, a contradiction.

Let $1\neq\gamma\in\Gamma_{2,v}$. Then $v$ does not lie on the mirror $M(\gamma)$, as otherwise any point on 
$e$ close enough to $v$ would be stabilised by $\gamma$. As this holds for all $\gamma\neq 1$ in $\Gamma_{2,v}$, we conclude that $\Gamma_{2,v}\cong C_p$. 
 Therefore, the vertex  $v$  has to $M(\gamma)$ the positive distance $v(\zeta_p-1)$ \cite[Lemma 3]{Her82}. 
 Analogously, $\Gamma_{2,w}\langle\delta\rangle\cong C_p$ and the distance between $w$ and $M(\delta)$
 is also $v(\zeta_p-1)$.
 
 Taking $e\subseteq\mathscr{T}_{\Gamma_1}^*$,
we  have $M(\gamma)\cup M(\delta)\subseteq\mathscr{T}_{\Gamma_1}^*$. From this, all three assertions follow.
\end{proof}

\begin{prop} \label{*tree4ptcycpp'}
Let $\Gamma\cong C_{pm}*C_{n}$ with $(n,p)=1$ be a discrete subgroup of $\PGL_2(K)$. Then
$T_\Gamma^*$ is as in Figure \ref{generalsegment} with $\dist(w,y)=0$. 
\end{prop}

\begin{proof}
The proof is similar to that of Proposition \ref{*tree4ptcycpp}.
\end{proof}

\begin{Rem}
The figure in \cite[Fig.\ 1]{BradMZ2005} is slightly erroneous. It should have
two segments stabilised by $C_{p^n}$
$$
\xymatrix@R=.5pt{C_{p^n}&C_{p^n}\\
*\txt{$\bullet$}\ar@{-}[r]^{C_{p^n}}
&*\txt{$\bullet$}\\
&
}
$$
which are not contained in the two mirrors. This can be seen by
 setting  $a=b=1$, and $r=s=n$ in Figure \ref{generalsegment}.
\end{Rem}

\subsubsection{Some examples.} 
\begin{figure}[h]
\centering
\includegraphics[scale=.27]{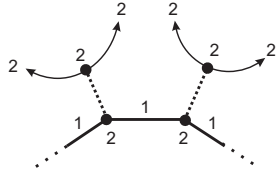}
\caption{The tree $\mathscr{T}_\Gamma^*$ for $\Gamma= C_2*C_2$ and $p=2$.
}
\label{c2c2-p=2}
\end{figure}

\begin{figure}[h]
\centering
\includegraphics[scale=.27]{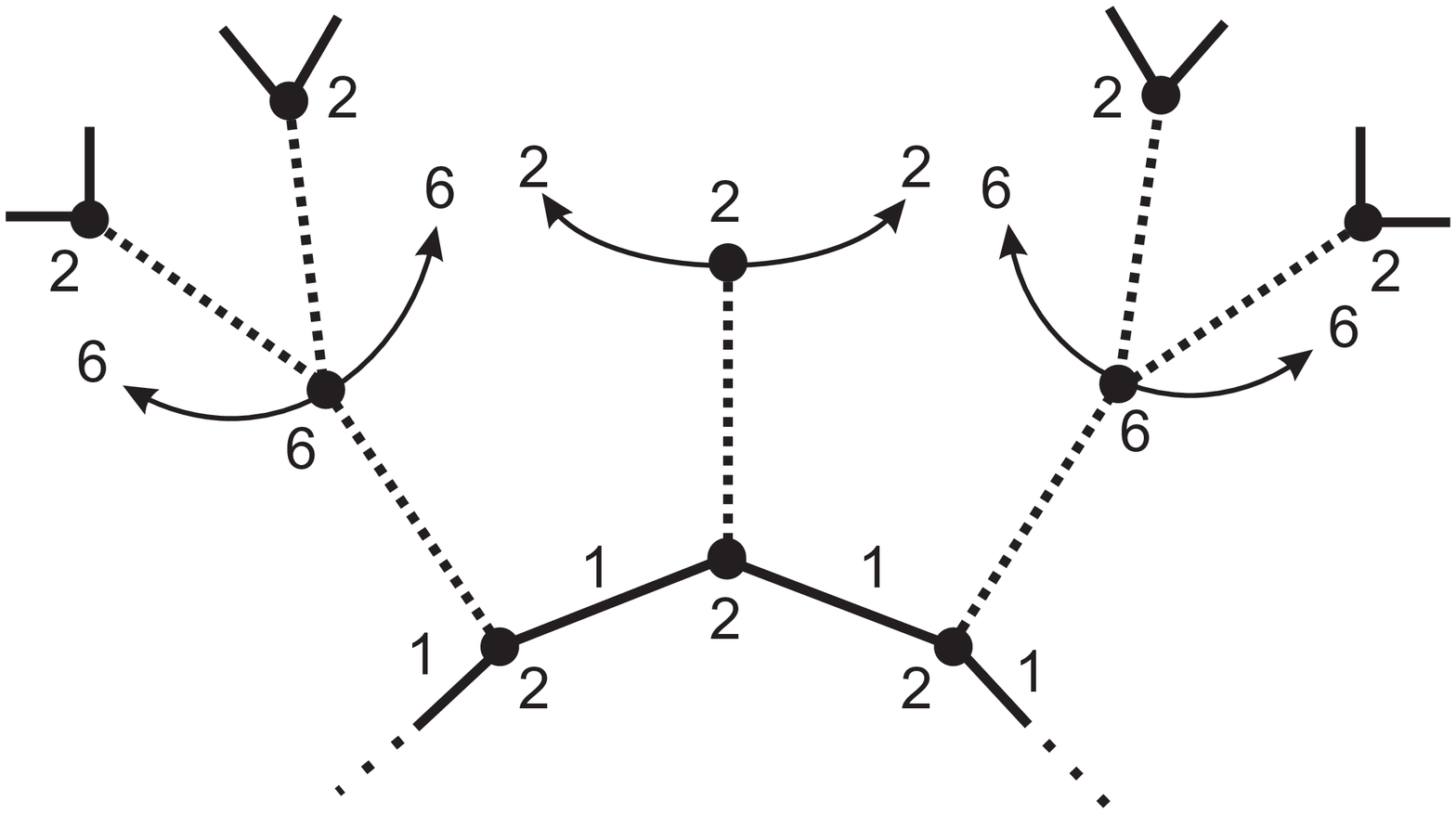}
\caption{The tree $\mathscr{T}_\Gamma^*$ for $\Gamma= C_2*C_6$ and $p=2$.
}
\label{c2c6-p=2}
\end{figure}

\begin{figure}[h]
\centering
\includegraphics[scale=.27]{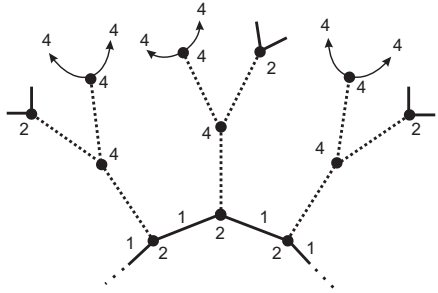}
\caption{The tree $\mathscr{T}_\Gamma^*$ for $\Gamma= C_4*C_4$ and $p=2$.
}
\label{c4c4-p=2}
\end{figure}

Figures \ref{c2c2-p=2}, \ref{c2c6-p=2} and \ref{c4c4-p=2} show portions of some $\mathscr{T}_\Gamma^*$ for $p=2$,
where $\dist(x,v)\neq 0$ (notation as in Figure \ref{generalsegment}).
We remark, however, that the most beautiful $*$-trees are those for the finite groups, when $p=2,3,5$
 (to appear in \cite{BKV01}),  some of which are illustrated already in \cite{CK2005}.

The figures are depicted in such a way that 
\begin{itemize}
\item a curved line with arrow heads represents the mirror
of a transformation whose order is written at both ends;
\item an unbroken line segment denotes an edge lying on a  geodesic line
on which a hyperbolic element acts through translation;
\item a dotted line segment means a   non-trivially stabilised edge which does not lie on a mirror
(the order of whose stabiliser is the lower of the two numbers at its extremities);
\item a number is the order of the stabiliser of the corresponding vertex or edge.
\end{itemize}

\subsubsection{ The positions of ends.}

Let $a=(a_0:a_1)$, $b=(b_0:b_1)$, $c=(c_0:c_1)$ and $d=(d_0:d_1)$ be four pairwise distinct $K$-rational points of $\mathbb{P}^1_K$. 
The arrangement of two straight lines $(a,b)$ and $(c,d)$ in  $\mathscr{T}_K$ can be calculated using
 the crossratio
$$
R(a,b;c,d)=\frac{(a_1c_0-a_0c_1)(b_1d_0-b_0d_1)}{(a_0b_1-a_1b_0)(c_0d_1-c_1d_0)}.
$$

\begin{prop} \label{crossratios}
Let $a,b,c,d\in\mathbb{P}^1_K$ be as above. Then:
\begin{enumerate}
\item If $|v(R(a,b;c,d))|=|v(R(b,a;c,d))|=0$, then $(a,b)$ and $(c,d)$ intersect at exactly one vertex. 
\item If $|v(R(a,b;c,d))|=|v(R(b,a;c,d))|\neq 0$, then $(a,b)$ and $(c,d)$ are disjoint with the distance $|v(R(a,b;c,d))|$. 
\item If $|v(R(a,b;c,d))|\neq|v(R(b,a;c,d))|$, then the intersection of $(a,b)$ and $(c,d)$ is the path $[v(a,b,c),v(b,c,d)]$ of length $$\max\{|v(R(a,b;c,d))|,|v(R(b,a;c,d))|\}.$$
\end{enumerate}
\end{prop}

Here, $v(a,b,c)$ denotes the unique vertex in $\mathscr{T}_K$ determined by the points $a,b,c\in\mathbb{P}^1_K$ viewed as ends in $\mathscr{T}_K$.

\begin{proof}
\cite[Proposition 3.5.1]{KatoJAG2005}.
\end{proof}

\begin{Def}
Let  $\zeta$ be a primitive $m$-th  root of unity. Then 
 $\epsilon_m$ and $\alpha_p(m,n)$
denote the numbers
\begin{align*}
\epsilon_m&:=\begin{cases}1,&\text{if $p\mid m$} \\0,&\text{otherwise}\end{cases}
\\[1mm]
\alpha_p(m,n)&:=\absolute{1-\zeta_p}^{\epsilon_m+\epsilon_n}.
\end{align*}
\end{Def}

$T_\Gamma^*$ has four ends which can be taken as $0,\infty$ going out of $x$ and $1,\lambda$ emanating from $y$ with $\absolute{\lambda}=1$.

\begin{thm} \label{branchpos}
Let $\Gamma\cong C_m*C_n$ be a discrete subgroup of $\PGL_2(K)$ and $0,\infty;1,\lambda$  the ends of $T_\Gamma^*$ as above. Then:
$$
\absolute{\lambda-1}<\alpha_p(m,n),
$$
and, conversely, for all such $\lambda\in K$ there is an embedding $C_m*C_n\to\PGL_2(K)$ as a discrete subgroup having such
a $*$-tree.
\end{thm}

\begin{proof}
Let $\Gamma$ be as stated.
Then we have, by Propositions \ref{*tree4ptcycp'p'}, \ref{*tree4ptcycpp} and \ref{*tree4ptcycpp'},
\begin{align*}
\dist(x,y)&=\dist(x,v)+\dist(v,w)+\dist(w,y)\\
&>
	\begin{cases}
	2\cdot v(\zeta_p-1), &\text{if $p\mid m$ and $p\mid n$,}\\
	v(\zeta_p-1),&\text{if $p\mid m$ and $(p,n)=1$,}\\
	0,&\text{otherwise},
	\end{cases}
\end{align*}
as $d(v,w)$ is strictly positive. Thus, by Proposition \ref{crossratios}, it holds true that
$$
\absolute{v(\lambda-1)}=\dist(x,y),
$$
as $\absolute{v(R(0,\infty;1,\lambda))}=\absolute{v(R(\infty,0;1,\lambda))}=\absolute{v(\lambda-1)}\neq 0$.
From this, the assertion follows.

For the converse implication, 
%
%
one has to check that the $*$-tree from Figure \ref{generalsegment}
is realisable for any value of $\dist(v,w)>0$ in $\absolute{v(K^\times)}$. 
This easy task is left to the reader. 
\end{proof}


\begin{exa}\rm
For $\Gamma\cong C_2*C_2$, we obtain the realisability for $\mathscr{T}^*_\Gamma$ if and only if
$\absolute{\lambda-1}<\absolute{2}^2$. In this case, the projective line is covered by a Tate elliptic curve.
By the formula for the $j$-invariant of elliptic curves \cite[Chapter IV.4]{Hartshorne},
$$
j=2^8\frac{(\lambda^2-\lambda+1)^3}{\lambda^2(\lambda-1)^2},
$$
this is equivalent to 
$
\absolute{j}>\absolute{2}^4.
$
\end{exa}

\subsection{Cyclic Mumford covers}

The only cyclic covers $X\stackrel{m}{\longrightarrow}\mathbb{P}^1$ allowing $X$ to be a Mumford curve are
known to be  those
corresponding to an equation of the form
\begin{align}\label{HMeq}
y^m=\prod\limits_{i=1}^r(x-\lambda_{i1})^{a_i}(x-\lambda_{i2})^{m-a_i}.
\end{align}
The branch points of the cover are the zeros of the polynomial on the right hand side.
After some projective $K$-linear transformation, the first four terms can be taken as
$x^{a_1},1,(x-1)^{a_2},(x-\lambda)^{m-a_2}$,
corresponding to the branch points $0,\infty;1,\lambda$. We will call a cover whose equation is
of the form (\ref{HMeq}), a {\em cover of  HM-type}.

\begin{dfn}
By an {\em $m$-cover}  of type $(e_1,\dots,e_r)$ we mean a cyclic cover $\phi\colon X\to\mathbb{P}^1$ of degree $m$ of HM-type ramified above the points
$(\lambda_{11},\lambda_{12};\dots;\lambda_{r1},\lambda_{r2})$ with 
$(\lambda_{11},\lambda_{12};\lambda_{21},\lambda_{22})=(0,\infty;1,\lambda)$
and $\absolute{\lambda}=1$, and such that the ramification index above each $\lambda_{ij}$
is $e_i>0$. 

A cyclic cover $X\to \mathbb{P}^1$ is called a {\em Mumford cover} or {\em of Mumford type}, if
$X$ is a Mumford curve.
\end{dfn}

\begin{dfn}
The statement
$$
\text{\em The bound holds for $m$.}
$$
means saying that an $m$-cover of type $(d,e)$ is a Mumford cover if and only if $\absolute{\lambda-1}<\alpha_p(m,n)$.
\end{dfn}

\begin{thm}
The bound holds for $m$.
\end{thm}

\begin{proof}
This, of course, is an immediate consequence of Theorem \ref{branchpos}. However, the statement will
be proven again in Section \ref{bound4m} by different methods.
\end{proof}

\subsection{Free products of cyclic groups}

\begin{lem} \label{lem-conjugfreeprod}
Let $G$ be a free product of finitely many groups $G_1,\dots, G_r$. Then each non-trivial element $s$ of finite order lies in
exactly one conjugate of  one of the factors $G_i$.
\end{lem}

\begin{proof}
By \cite[IV.1.6]{LS1977}, $s$ lies only in conjugates of some of the $G_i$. Assume therefore that $s\in G_1$. Then the equation $s=g^{-1}s_ig$ with $s_i\in G_i$ and $i\neq 1$ is easily seen to lead to a contradiction.
\end{proof}

Let $N=\langle s_0\rangle*\dots*\langle s_m\rangle\subseteq\PGL_2(K)$ be the $m$-fold free product of the 
cyclic group $C_n$ acting discontinuously on $\mathbb{P}^1$. By the universal property of free products, there is  a unique homomorphism $\phi\colon N\to C_n$ such that for each $i=0,\dots,m$ the diagram
$$
\xymatrix{
&\langle s_i \rangle \ar[dl] \ar[dr]^\cong&\\
N\ar[rr]_{\phi}&&C_n
}
$$
is commutative. This homomorphism $\phi$ depends on the choices 
of the isomorphisms $\langle s_i\rangle\to C_n$. Here, all $s_i$ are
supposed to be mapped to the same generator of $C_n$. In Section \ref{bound4m},
we will consider also other choices.

\medskip
Let $\Gamma:=\ker\phi\subseteq N$. It  is easily seen to be of finite index $n$ in $N$.

\begin{Prop}	\label{thm-Reidemeistergen}
The  group\ $\Gamma$ is free of rank $m(n-1)$ and freely generated by
$$
s_0^j s_i s_0^{-j-1},\qquad j=1,\dots,n-1,\quad i=1,\dots,m.
$$
\end{Prop}

\begin{proof}
By \cite[\S1.3 and \S2.(4)]{Reidemeister1927}, $\Gamma$ is generated by the asserted elements. This generating system cannot
be shortened, as the genus $g$ of the Mumford curve $X=\Omega_N/\Gamma$ can be calculated by the Riemman-Hurwitz formula for the cyclic cover 
$$
X\to\mathbb{P}^1
$$
of degree $n$ 
 totally ramified in the  $\Gamma$-orbits of the points in $\Omega_N$ fixed by some $s_i$.
For this, we must check that the $s_i$ have regular fixed points: this follows from 
\cite[Satz 6]{Herrlich1980}, as, by Lemma \ref{lem-conjugfreeprod},
$s_i$ fixes precisely one vertex of $\mathscr{T}_N$. 
So, we have $2m+2$ ramification points, and
$$
2g-2=-2n+(2m+2)(n-1),
$$  
which is equivalent to 
$$
g=m(n-1).
$$

\smallskip
Any non-trivial element  $\gamma\in\Gamma$ of finite order is conjugated to an element of some 
$\langle s_i\rangle$. It follows that,  in any representation of $\gamma$
as a word in the generators $s_0^j s_i s_0^{-j-1}$ and their inverses, the sum of the exponents cannot be  zero---a contradiction. Therefore, $\Gamma$ is torsion-free. As $\Gamma$ is the fundamental group of
a tree of groups, it follows that $\Gamma$ is free. 
\end{proof}

\section{Many-point Mumford covers}

\begin{Lemma}	\label{discreteedges}
Let $N\subseteq \PGL_2(K)$ be a free tree product of finite cyclic groups $C_1,\dots, C_r$. Then $N$
is discrete if and only if each free amalgam $C_i*C_j\subseteq N$ of two  neighbouring factors of $N$ 
 is discrete. 
\end{Lemma}

\begin{proof}
If $N$ is discrete, then so is the  subgroup $C_i*C_j$.

We prove the converse by induction on $r$. If $r=2$, then the  statement clearly holds true.
Let, for $r>2$, $N=N'*C$ with $N'$ a free tree product with $r-1$ factors and $C=C_i$ for some $i$
between $1$ and $r$. By the induction hypothesis, $N'$ is discrete. Also $C':=C*C_j$  is discrete, 
where $C_j$ is the unique factor of $N'$ neighbouring to $C$.  
Clearly, 
$
T:=T_{N'}^*\cup T_{C'}^*\subseteq\mathscr{T}_K
$ 
is a tree, and 
$$
\mathscr{T}:=\bigcup\limits_{\gamma\in N}\gamma\, T\subseteq \mathscr{T}_K
$$ 
is a tree
upon which $N$ acts with finite vertex stabilisers: 
$$
N_v\cong\begin{cases}
	N'_v,&v\in\gamma\, T^*_{N'}\\
	C'_v,&v\in\gamma\,T^*_{C'}
	\end{cases}		
\quad\text{for some $\gamma\in N$.}
$$
By \cite[Lemma 4.4.1(2)]{KatoJAG2005}, it follows that $N$ is discrete.
\end{proof}

\begin{rem}\rm
In fact, with the notations from the proof of Lemma \ref{discreteedges}, it holds true that
$$
\mathscr{T}=\mathscr{T}^*_N.
$$
This is due to the fact that $N$ is generated by the stabilisers $N_v$, where $v$ runs through all vertices
of $T$ (cf.\ the "if`` part in "$g=0$`` of the proof   of Theorem II in \cite[\S7.]{KatoJAG2005}).
\end{rem}

Lemma \ref{discreteedges} allows us to prove a geometric criterion for arbitrary $m$-covers to
be Mumford covers. For this, denote by $\branchlocus(\phi)$ the branch locus of an $m$-cover.

\begin{thm} \label{geocond}
An $m$-cover of type $(e_1,\dots,e_r)$ is a Mumford cover, if and only if, after a suitable re-ordering of
the pairs $(\lambda_{ij},e_i)$,  there is an affinoid covering $\mathfrak{U}=\{U_1,\dots,U_r\}$ of $\mathbb{P}^1$
such that 
\begin{enumerate}
\item  $U_i\cap U_j$ is either empty or an annulus of thickness $\alpha_p(e_i,e_j)$, if $i\neq j$,
\item for all $i=1,\dots,r$ holds true: $U_i\cap\branchlocus(\phi)=\{\lambda_{i1},\lambda_{i2}\}$.
\end{enumerate}
\end{thm}

\begin{proof}
If an $m$-cover  $\phi\colon X\to\mathbb{P}^1$ is a Mumford cover, then $\mathbb{P}^1$ can be uniformised by a free tree product $N$ of
cyclic groups,
i.e.\ $\phi$ is part of a commutative diagram
$$
\xymatrix{
\Omega\ar[dr]_N\ar[r]&X\ar[d]^\phi\\
&\mathbb{P}^1
}
$$ 
where $X$ is a Mumford curve, and 
the space $\Omega\subseteq\mathbb{P}^1$ is the
complement of the closure of the set of limit points for
the action of $N$ given by some
discrete faithful representation $\tau\colon N\to\PGL_2(K)$. 

The tree $\mathscr{T}^*_N$, embedded in $\mathscr{T}_K$ via $\tau$, allows the extraction of the cover $\mathfrak{U}$: the stars around the vertices
whose stabilisers are maximal yield  discs $U_i$ with $\deg(v)-1$ "holes``, and the paths between
any two nearest such vertices correspond to annuli whose thickness was calculated in the proof of Theorem \ref{branchpos}
as $\alpha_p(e_i,e_j)$. 
 
Let, conversely, $\phi$ be an $m$-cover satisfying the conditions (1) and (2). Taking two intersecting
$U_i, U_j\in\mathfrak{U}$, we can construct a four-point $m$-cover by setting the ramification indices of all
branchpoints  outside $U_i\cup U_j$ to one. This is a Mumford cover, by Theorem \ref{branchpos}.
Doing this for all intersecting pairs of affinoids from $\mathfrak{U}$, we obtain a free amalgamated  product
which is discrete, by Lemma \ref{discreteedges}. 
\end{proof}

\begin{cor}\label{cor-geocond}
A cyclic cover $\phi\colon X\to\mathbb{P}^1$ is of Mumford type, if and only if there is some  $\alpha\in\PGL_2(K)$ such that 
$\alpha\circ\phi$  is
an $m$-cover satifying conditions (1) and (2) of Theorem \ref{geocond}.
\end{cor}

\begin{proof}
This follows immediately from the fact that the cross-ratio of any four points in $\mathbb{P}^1$ is invariant under
projective linear transformations.
\end{proof}

\begin{rem}\rm
Theorem \ref{geocond} generalises the characterisation 
from \cite{Ulrich1981} of hyperelliptic Mumford curves among $2$-covers, 
 proven in the case of residue characteristic unequal $2$ and
by entirely different methods.
This geometric condition is used by Frank Herrlich for constructing a moduli space of hyperelliptic Mumford curves \cite{Herrlich1986}.
\end{rem}

\section{The bound for four-point covers again} \label{bound4m}

In the following, we will calculate the bound for four-point cyclic covers
 by giving explicit faithful representations 
$\tau\colon N=C_m*C_n\stackrel{~}{\longrightarrow}\langle s,t\rangle\subseteq\PGL_2(K)$. 
In the discrete case, the 
fixed points of $s$ and $t$ correspond then to four ``upstairs" ramification  points  of
the cover $\xymatrix{\Omega_{\tau(N)}\ar[r]^{\tau(N)}&\mathbb{P}^1}$ which we may and will assume to be
$0,\infty;1,\lambda$ with $\absolute{\lambda}=1$. 

\smallskip
The following Lemma shows that this approach is  indeed  legitimate, albeit  indirect.
 
\begin{lem} \label{firstorderbd}
Assume that the branch locus of $\phi$ is $0,\infty,1,\lambda'$ with $\absolute{\lambda'}=1$.
Then it holds true that
$$
\absolute{\lambda-1}<\alpha_p(m,n)\iff\absolute{\lambda'-1}<\alpha_p(m,n).
$$
\end{lem}

\begin{proof}
This follows from the fact that any section $T_N^*\to\mathscr{T}_N^*$ is isometric.  
\end{proof}

\begin{rem}
In fact the first approach from Section \ref{mumfdiscr} was indirect in the same manner as the 
approach in this section, as we calculated $T_N^*$ {\em within} $\mathscr{T}^*_N$.

The Kummer equations $y^m=x^a(x-1)^b(x-\lambda)^c$ to follow are to be understood 
modulo Lemma \ref{firstorderbd}.
of the corresponding covers. In fact, the precise correspondence between discrete faithful representations
of $N$ and Kummer equations is still not settled. 
\end{rem}

\subsection{Galois covers of prime degree} \label{explicitprime}


Let $X\to\mathbb{P}^1$ be a cyclic cover of prime degree $q$ totally ramified above exactly four points. By projective linear transformation, we may assume that the branch locus
of the cover consists of the points $0$, $1$, $\infty$ and $\lambda$, where $\absolute{\lambda}=1$. The aim
of this section is to redetermine explicitly the conditions on  $\lambda$ for which $X$ can be a Mumford curve by using Ford's isometric circles\footnote{For the notion of {\em isometric circles} and their properties, cf.\ \cite[Ch.\ I, \S 11]{Ford1927}.}.


Let $N=N_{q,q}=\langle s\rangle * \langle t\rangle\subseteq\PGL_2(K)$ be the free product of two copies of the cyclic group $C_q$,
where $s$ is given by the  matrix
$$
s=\begin{pmatrix}\zeta&0\\0&1\end{pmatrix},
$$
where $\zeta$ is a primitive $q$-th root of unity, and $t$ is obtained from $s$ by conjugation with
$$
\phi=\begin{pmatrix}\lambda&1\\1&1\end{pmatrix}. 
$$
The latter means that
$$
t=\phi s\phi^{-1}=\frac{1}{\lambda-1}\begin{pmatrix}\lambda\zeta-1&\lambda(1-\zeta)\\ \zeta-1&\lambda-\zeta\end{pmatrix}.
$$
The elliptic transformation $s$ has the fixed points $0$ and $\infty$, whereas $t$ has the fixed points $1$ and
$\lambda$. 

For further reference, we also give the matrix $t^{-1}$:
$$
t^{-1}=\phi s^{-1}\phi^{-1}=\frac{1}{\lambda-1}
\begin{pmatrix}\lambda\zeta^{-1}-1&\lambda(1-\zeta^{-1})\\\zeta^{-1}-1&\lambda-\zeta^{-1}\end{pmatrix}.
$$

\begin{Lemma} \label{schottkygen-q}
The normal free subgroups $\Gamma$ of $N$ of index $q$ are all of rank $q-1$ and given as $\Gamma=\Gamma_f=\ker\phi_f$ ($f=1,\dots,q-1$), 
where each $\phi_f$ is the map
$$
\phi_f\colon N\to C_q,\quad s\mapsto\zeta,\;t\mapsto \zeta^f.
$$
\end{Lemma}

\begin{proof}
The $\Gamma_f=\ker\phi_f$ are clearly normal
and, by Proposition \ref{thm-Reidemeistergen}, these groups are free of rank $q-1$.
These are in fact all normal subgroups of index $q$, as every group homomorphism
$\phi\colon N\to C_q$ factorises through the abelian group $N^\abelian=\langle s\rangle\times\langle t\rangle$:
$$
\xymatrix{
N\ar[rr]^\phi \ar[dr] && C_q\\
&N^\abelian \ar[ur]_\psi&
}
$$
and the map $\psi$ can be made into the form 
$$
\psi_f\colon N^\abelian\to C_q,\quad s\mapsto \zeta,\;t\mapsto \zeta^f
$$ 
via an automorphism of $C_q$.
\end{proof}

\begin{Thm}	\label{thm-primexplicit}
The equation
$$
y^q=x(x-1)^a(x-\lambda')^b,
$$
where $1\le a,b<q$ and $\absolute{\lambda'}=1$, defines a covering of the projective line by a Mumford
curve whose topological fundamental group 
is $\Gamma_f$  as in Lemma \ref{schottkygen-q}, if and only if 
$$
b=q-a\quad \text{and}\quad \absolute{\lambda'-1}<\alpha_p(q,q)
=
\begin{cases}
\absolute{1-\zeta_p}^2,&p=q\\
1,&\text{otherwise.}
\end{cases}
$$
Here, the number $f$  is such that $af\equiv 1\mod q$.
\end{Thm}

\begin{proof}
The condition $b=q-a$ and $af\equiv 1\mod q$  on the exponents was found by van Steen   using theta functions \cite[Proposition 3.2]{vanSteen1982}. 

\smallskip
The generators for $\Gamma_f$ from Proposition \ref{thm-Reidemeistergen} are
$$
\gamma_{if}=s^its^{-f-i}=
\frac{1}{\lambda-1}\begin{pmatrix}(\lambda\zeta-1)\zeta^{-f}&\lambda\zeta^i(1-\zeta)\\
							\zeta^{-i-f}(\zeta-1)&\lambda-\zeta\end{pmatrix},
$$
where $i=1,\dots,q-1$. An automorphism $\gamma$ of $\mathbb{P}^1$ is hyperbolic if and only if
$\frac{\absolute{\Trace\gamma}^2}{\absolute{\det\gamma}}>1$. Now, 
$$
\Trace\gamma_{if}=\frac{(1+\zeta^{1-f})\lambda-(\zeta+\zeta^{-f})}{\lambda-1},\qquad \det\gamma_{if}=\zeta^{1-f},
$$ 
therefore,
\begin{align*}
\frac{\absolute{\Trace\gamma_{if}}^2}{\absolute{\det\gamma_{if}}}>1&\\
\iff
\absolute{\lambda-1} <\; \absolute{(1+\zeta^{1-f})\lambda-&(\zeta+\zeta^{-f})}\\
=\absolute{1+\zeta^{1-f}-&(\zeta+\zeta^{-f})}\leq 1,
\end{align*}
where the equality holds, because the difference of the two corresponding terms has norm $\absolute{\lambda-1}$ or less. It follows that, in the case 
$$
\absolute{\lambda-1}=\absolute{1+\zeta^{1-f}-(\zeta+\zeta^{-f})}=\absolute{1-\zeta}\absolute{1-\zeta^{-f}}
=\alpha_p(q,q),
$$
the group $\Gamma_f$ is not discontinuous and therefore does not give rise to a Mumford curve.

\smallskip
Let us assume that
$
\absolute{\lambda-1}<\absolute{1-\zeta}^2.
$ 
The isometric circles for $\gamma_{if}$ and 
$$
\gamma_{if}^{-1}=s^{f+i}ts^{-i}=\frac{1}{\lambda-1}
\begin{pmatrix}
(\lambda\zeta^{-1}-1)\zeta^f&\lambda\zeta^{f+i}(1-\zeta^{-1})\\(\zeta^{-1}-1)\zeta^{-i}&\lambda-\zeta^{-1}\end{pmatrix}
$$
are 
\begin{align*}
I_{\gamma_{if}}&=\left\{z\in\mathbb{P}^1:\left|z-\frac{\zeta-\lambda}{\zeta-1}\zeta^{i+f}\right|
<\frac{\absolute{\lambda-1}}{\absolute{\zeta-1}}\right\},
\\
I_{\gamma_{if}^{-1}}&=\left\{z\in\mathbb{P}^1:\left|z-\frac{\zeta^{-1}-\lambda}{1-\zeta}\zeta^{1-i}\right|<\frac{\absolute{\lambda-1}}{\absolute{1-\zeta}}\right\}.
\end{align*}
One then sees  that 
\begin{align*} \label{nointersect}
I_{\gamma_{if}}^+\cap I_{\gamma_{if}^{-1}}^+=I_{\gamma_{if}}^+\cap I_{\gamma_{jf}^{}}^+
=I_{\gamma_{if}^{-1}}^+\cap I_{\gamma_{jf}^{-1}}^+=I_{\gamma_{if}}^+\cap I_{\gamma_{jf}^{-1}}^+
=\emptyset
\end{align*}
for all $i,j=1,\dots,q-1$.
Therefore,  the complement of the union of these open disks is a good fundamental domain for $\Gamma_f$ in the sense of \cite[(4.1.3)]{GvP1980}.
\end{proof}

\begin{exa}\rm
Specialising the calculations for $q=2$, one obtains again
$$
\absolute{\lambda-1}<\absolute{2\lambda+2}\leq 1\implies \absolute{\lambda-1}<\absolute{2}\absolute{\lambda+1}=\absolute{2}^2,
$$
because indeed $\absolute{\lambda+1}=\absolute{1+1}$, due to $\absolute{\lambda-1}<1$.
\end{exa}

\subsection{Totally ramified four-point covers}

We assume that $\phi\colon X\to \mathbb{P}^1$ is of degree $m$ and totally ramified above the four branch points.
Let $N_{m,m}=\langle s\rangle*\langle t\rangle$ with $s$ and $t$ of order $m$. 

\begin{Thm} \label{descent1}
Let $q$ be a prime dividing $m$. Then the  bound holds for $m$ if  
it holds for  $m'=\frac{m}{q}$.  
\end{Thm}

\begin{proof}
Assume that the bound holds for $m'$. We know already that it holds for $q$. 
Therefore, the diagram
with exact rows and columns (and $\Gamma_{q,q}=\Gamma_f$ as in the preceding subsection)
$$
\xymatrix{
&1\ar[d]&1\ar[d]&1\ar[d]&\\
1\ar[r]&\Gamma_{q,q}\ar[r]\ar[d]&N_{q,q}\ar[d]\ar[r]&C_q\ar[d]\ar[r]&1\\
1\ar[r]&\Gamma_{m,m}\ar[r]\ar[d]&N_{m,m}\ar[d]\ar[r]&C_m\ar[d]\ar[r]&1\\
1\ar[r]&\Gamma_{m',m'}\ar[r]\ar[d]&N_{m',m'}\ar[d]\ar[r]&C_{m'}\ar[d]\ar[r]&1\\
&1&1&1&\\
}
$$ 
yields generators for $\Gamma_{m,m}$ which can be examined by the method of isometric circles.
Indeed, 
$$
\Gamma_{m,m}=\langle \Gamma_{q,q}, \zeta_q^i\Gamma_{m',m'}\zeta_q^{-i}\mid i=0,\dots,q-1\rangle,
$$
where $\zeta_q$ is a generator of $C_q$.
By assumption, both $\Gamma_{q,q,}$ and $\Gamma_{m',m'}$ are free of rank $q-1$ and $m'-1$, respectively.
As the right and middle columns are split, also the left column splits. Therefore,
$\Gamma_{m,m}$ 
is free of
 rank 
 $$
g_{m,m}=(q-1)+(m'-1)q=m-1.
$$
The generators obtained in this way from generators of $\Gamma_{q,q}$ and $\Gamma_{m',m'}$ are hyperbolic if and only if the latter
are both Schottky groups, 
which   is equivalent to 
$$
\absolute{\lambda-1}<\min\{\alpha_p(q,q),\alpha_p(m',m')\},
$$ 
by assumption.
But then one calculates that the isometric circles of any pairs of different generators of $\Gamma_{m,m}$ and their
inverses do not intersect. Thus the bound holds for $m$.
\end{proof}

\begin{cor} \label{bound4totram}
The bound holds for $m$, if the $m$-cover is totally ramified.
\end{cor}

\begin{proof}
This follows by an iterative application of Theorem \ref{descent1}.
\end{proof}

\subsection{Four point covers with arbitrary ramification}

Let $X\to\mathbb{P}^1$ be a cover of degree $n$ ramified above the points $0,1,\infty,\lambda'$ with $\absolute{\lambda'}=1$
given by the pair $(N,\Gamma)$ with $N=\langle s\rangle*\langle t\rangle$ and a free normal subgroup $\Gamma$ realised as the kernel of a surjection $N\to C_n$. Let the orders of $s$ and $t$ be $d$ and $e$.

\subsubsection{The case $e\mid d$}

In this case, $n=d$, as otherwise $X$ would not be connected. Let $\zeta$ be a primitive $n$-th root of unity, and $f:=\frac{n}{e}$. 
As before, consider the maps
$$
\phi_k\colon N\to C_n,\quad s\mapsto\zeta,\;t\mapsto\left(\zeta^f\right)^k,\qquad (k,e)=1.
$$
The same method by Reidemeister as before yields generators for $\Gamma_k=\ker\phi_k$
$$
B_k:=\{\gamma_{ijk}=s^{ik}\gamma_{jk} s^{-ik}\mid i=1,\dots,f,\;j=1,\dots,e-1\},
$$
where
$$
\gamma_{jk}=\left(s^{fk}\right)^j\:t\left(s^{fk}\right)^{-j-1}.
$$

\begin{Thm}
The bound holds for $m$-covers of type $(d,m)$.
\end{Thm}

\begin{proof}
Let $m=d\ell$ and
consider the commutative diagram 
$$
\xymatrix{
&X\ar[d]^{C_d}\\
\Omega\ar[r]^{\ast C_d}\ar[dr]_{N_{d,m}}\ar[ur]^{\Gamma_{d,m}}&\mathbb{P}^1\ar[d]^{C_\ell}\\
&\mathbb{P}^1
} 
$$
The vertical maps   $\psi\colon X\to\mathbb{P}^1$
and $\phi\colon \mathbb{P}^1\to\mathbb{P}^1$
are cyclic, with a Mumford curve $X$, and
$\phi$ is ramified above $\mathset{0,\infty}$.  
The branch locus of the horizontal map 
is 
$f^{-1}(\mathset{0,\infty,1,\lambda})$ which coincides with the branch locus of 
$\psi$ and is of cardinality $2\ell+2$.
By looking at
the corresponding $*$-trees, 
we see that $T^*_{\ast C_d}$ has exactly $\ell+1$ vertices
stabilised by $C_d$: more precisely,
from one  vertex $v$ on one mirror,
there are  $\ell$ paths to the other mirrors, and the pairwise intersection
of thes paths is $v$ (in other words, $T_{\ast C_d}$ is star-shaped
with centre $v$).
Hence, $\ast C_d$ is a free tree product of $\ell+1$ copies of $C_d$.

Now, the top triangle with the cyclic cover $\psi$ yields
that $\Gamma_{d,m}$ is isomorphic to a free product of $\ell$
copies of $\Gamma_{d,d}$,
which is part of an exact sequence
$$
\xymatrix{
1\ar[r]&\Gamma_{d,d}\ar[r]&C_d*C_d\ar[r]&C_d\ar[r]&1
}
$$  
However, from the proof of Theorem \ref{descent1} we know that $\Gamma_{d,d}$
is free of rank $d-1$.
For $C_\ell=\langle\zeta\rangle$, this implies that
$$
\Gamma_{d,m}=\langle\zeta^i\Gamma_{d,d}\zeta^{-i}\mid i=0,\dots,\ell-1\rangle
$$
  is free of rank
$$
g_{d,d\ell}=\ell(d-1).
$$
By Corollary \ref{bound4totram}, the bound holds for  
$\psi$. The section $T^*_{N_{d,m}}\to T^*_{\ast C_d}$
being isometric implies that the bound holds for $\phi\circ\psi$. 
\end{proof}

\subsubsection{The case $(d,e)=1$} 	\label{decoprime}

In the case that $(d,e)=1$, it follows that necessarily $n=ed$. Consider the maps
$$
\phi_{k\ell}\colon N\to C_n,\quad s\mapsto\zeta^{ek},\; t\mapsto\zeta^{d\ell},
$$
where $(k,e)=1$ and $(\ell,d)=1$. Let 
$\sigma:=s^e$, $\tau:=t^d$
and
$$
B_{k\ell}:=\{ \gamma_{ij}:=\sigma^{-i}\tau^{-j}\sigma^{i}\tau^{j}
\mid i=1,\dots,d-1,\;j=1,\dots,e-1\}.
$$

\begin{Prop}
$\Gamma_{k\ell}:=\ker\phi_{k\ell}$ is free of rank $(e-1)(d-1)$.
\end{Prop}

\begin{proof}
This follows from a similar Riemann-Hurwitz argument as in the proof of Proposition \ref{thm-Reidemeistergen}.
\end{proof}

Now assume that $(p,d)=1$.

\begin{Thm}	\label{thm-edcoprime}
The equation
$$
y^n=x^a(x-1)^b(x-\lambda')^{n-b},
$$
where  $1\le a< n$ is of order $e$ mod $n$, $1\le b< n$ of order $ d$ mod $n$,
$(d,e)=1$ and $\absolute{\lambda'}=1$ defines a Mumford curve covering $\mathbb{P}^1$,
if and only if $\absolute{\lambda'-1}<\alpha_p(1,e)$. 
\end{Thm}

\begin{proof}
From
\begin{align*}
\sigma^i\tau^j=&\frac{1}{\lambda-1}
	\begin{pmatrix}
	(\lambda\zeta^{dj}-1)\zeta^{ei}&\lambda(1-\zeta^{dj}\zeta^{ei})\\
	\zeta^{dj}-1 & \lambda-\zeta^{dj}
	\end{pmatrix},
\\[2mm]
\sigma^{-i}\tau^{-j}=&\frac{1}{\lambda-1}
	\begin{pmatrix}
	(\lambda\zeta^{-dj}-1)\zeta^{-ei}&\lambda(1-\zeta^{-dj}\zeta^{-ei})\\
	\zeta^{-dj}-1 & \lambda-\zeta^{-dj}
	\end{pmatrix}
\end{align*}
we calculate
\begin{align*}
\gamma_{ij}:=\sigma^{-i}\tau^{-j}\sigma^i\tau^j=
\frac{1}{(\lambda-1)^2}
\begin{pmatrix}
a_{ij}&b_{ij}\\c_{ij}&d_{ij}
\end{pmatrix}
\end{align*}
with
\begin{align*}
a_{ij}&=(\lambda\zeta^{dj}-1)(\lambda\zeta^{-dj}-1)-\lambda\zeta^{-ei}(\zeta^{dj}-1)(\zeta^{-dj}-1),
\\
b_{ij}&=\lambda(\lambda\zeta^{-dj}-1)(1-\zeta^{dj})+\lambda\zeta^{-ei}(1-\zeta^{-dj})(\lambda-\zeta^{dj}),
\\
c_{ij}&=(\zeta^{-dj}-1)(\lambda\zeta^{dj}-1)\zeta^{ei}+(\lambda-\zeta^{-dj})(\zeta^{dj}-1),
\\
d_{ij}&=\lambda\zeta^{ei}(\zeta^{-dj}-1)(1-\zeta^{dj})+(\lambda-\zeta^{-dj})(\lambda-\zeta^{dj}).
\end{align*}
By definition, $\det\gamma_{ij}=1$. The trace of $\gamma_{ij}$ is
$$
\Trace\gamma_{ij}=\frac{2(\lambda-\zeta^{-dj})
	(\lambda-\zeta^{dj})-\lambda(\zeta^{-ei}+\zeta^{ei})(1-\zeta^{dj})(1-\zeta^{-dj})}
{(\lambda-1)^2}.
$$
Thus the condition for hyperbolicity of $\gamma_{ij}$ is
$$
\absolute{\lambda-1}^2<\absolute{2(\lambda-\zeta^{-dj})
	(\lambda-\zeta^{dj})-\lambda(\zeta^{-ei}+\zeta^{ei})(1-\zeta^{dj})(1-\zeta^{-dj})}.
$$
Set
$$
\epsilon_{ij}:=(1-\zeta^{dj})(1-\zeta^{-dj})(1-\zeta^{ei})(1-\zeta^{-ei}),
$$
and notice that
$$
\zeta^{ei}+\zeta^{-ei}=2-(1-\zeta^{ei})(1-\zeta^{-ei}).
$$
Therefore,
the right hand side of the inequality equals
\begin{align*}
&
\absolute{2((\lambda-\zeta^{-dj})(\lambda-\zeta^{dj})-\lambda(1-\zeta^{dj})(1-\zeta^{-dj}))+\lambda\epsilon_{ij}}
\\
=&
\absolute{2(\lambda^2-\lambda(\zeta^{dj}+\zeta^{-dj})+1-2\lambda+\lambda(\zeta^{dj}+\zeta^{-dj}))+\lambda\epsilon_{ij}}
\\
=&\absolute{2(\lambda-1)^2+\lambda\epsilon_{ij}}
\\
=&\absolute{\lambda\epsilon_{ij}}=\absolute{\epsilon_{ij}},
\end{align*}
where the first equality in the last line holds true, because $\absolute{2(\lambda-1)^2}\leq\absolute{\lambda-1}^2$.

\smallskip
In a similar way we obtain
$$
\gamma_{ij}^{-1}=\tau^{-j}\sigma^{-i}\tau^j\sigma^i=\frac{1}{\lambda-1}
	\begin{pmatrix}
	a'_{ij}& b_{ij}'\\
	c'_{ij}&d'_{ij}
	\end{pmatrix},
$$
with
\begin{align*}
a'_{ij}&=(\lambda\zeta^{dj}-1)(\lambda\zeta^{-dj}-1)-\lambda\zeta^{ei}(\zeta^{dj}-1)(\zeta^{-dj}-1),
\\
b'_{ij}&=\lambda\zeta_{-ei}(\lambda\zeta^{-dj}-1)(1-\zeta^{dj})+\lambda(1-\zeta^{-dj})(\lambda-\zeta^{dj}),
\\
c'_{ij}&=(\zeta^{-dj}-1)(\lambda\zeta^{dj}-1)+\zeta^{ei}(\lambda-\zeta^{-dj})(\zeta^{dj}-1),
\\
d'_{ij}&=\lambda\zeta^{-ei}(\zeta^{-dj}-1)(1-\zeta^{dj})+(\lambda-\zeta^{-dj})(\lambda-\zeta^{dj}).
\end{align*}

The isometric circles are
\begin{align*}
I_{\gamma_{ij}}&=\left\{z\in\mathbb{P}^1:
		\left|z+\frac{d_{ij}}{c_{ij}}\right|<\frac{\absolute{\lambda-1}^2}{\absolute{c_{ij}}}
	\right\}
\\
I_{\gamma_{ij}^{-1}}&=\left\{z\in\mathbb{P}^1:
		\left|z+\frac{d'_{ij}}{c'_{ij}}\right|<\frac{\absolute{\lambda-1}^2}{\absolute{c'_{ij}}}
	\right\}.
\end{align*}
They do not intersect pairwise,
if and only if
$$
\absolute{\lambda-1}^2<\min\left\{\absolute{d_{ij}-d'_{ij}},\;\absolute{d_{ij}-d_{i'j'}},\;\absolute{d'_{ij}-d'_{i'j'}}
\right\},
$$
where $i,i'=1,\dots,d-1$ and $j,j'=1,\dots,e-1$ are such that the set to be minimised does not contain zero.

Rewrite $d_{ij}$ as
\begin{align*}
d_{ij}&=(\lambda-\zeta^{-dj})(\lambda-\zeta^{dj})-\lambda\zeta^{ei}(\zeta^{-dj}-1)(\zeta^{dj}-1)
\\
&=\lambda^2-\lambda(\zeta^{-dj}+\zeta^{dj})+1-\lambda\zeta^{ei}\left(2-(\zeta^{dj}+\zeta^{-dj})\right)-2\lambda+2\lambda
\\
&=(\lambda-1)^2+2\lambda(1-\zeta^{ei})+\lambda(\zeta^{dj}-\zeta^{dj})(\zeta^{ei}-1)
\\
&=(\lambda-1)^2+\lambda(1-\zeta^{ei})\left(2-(\zeta^{dj}+\zeta^{-dj})\right)
\\
&=(\lambda-1)^2+\lambda(1-\zeta^{ei})(1-\zeta^{dj})(1-\zeta^{-dj}),
\end{align*}
and, similarly, $d'_{ij}$ 
as
$$
d_{ij}'=(\lambda-1)^2+\lambda(1-\zeta^{-ei})(1-\zeta^{dj})(1-\zeta^{-dj}),
$$
and set $e=p^r\ell$, $(p,\ell)=1$. Then
the minimum is attained for $j=j'=p^{r-1}\ell$
and takes the value
$$
\absolute{d_{ij}-d_{i'j}}=\absolute{\lambda}\cdot\absolute{\zeta^{ei'}-\zeta^{ei}}\cdot\absolute{(1-\zeta^{dj})(1-\zeta^{-dj})}
=\absolute{1-\zeta_p}^2\,,
$$
since we assumed $(d,p)=1$.
\end{proof}

\subsubsection{The case $d\not|\;\,e$ and $e\not|\;\, d$}

In the case $d\not|\;\,e$ and $e\not|\;\, d$, we have
$$
n=\lcm(d,e)\quad\text{and}\quad \ell=\gcd(d,e). 
$$

\begin{Thm}
The equation 
$$
y^n=x^a(x-1)^b(x-\lambda')^{n-b}
$$
where  $1\le a< n$ is of order $e$ mod $n$, $1\le b< n$ of order $ d$ mod $n$,
and $\absolute{\lambda'}=1$ defines a Mumford curve covering $\mathbb{P}^1$,
if and only if 
$$
\absolute{\lambda'-1}<\alpha_p(d,e).
$$ 
\end{Thm}

\begin{proof}
 Let $d':=\frac{d}{\ell}$, $e':=\frac{e}{\ell}$, 
$m:=d'e'$ and consider the diagram
$$
\xymatrix{
&1\ar[d]&1\ar[d]&1\ar[d]&\\
1\ar[r]&\Gamma_{d',e'}\ar[r]\ar[d]&N_{e',d'}\ar[d]\ar[r]&C_m\ar[d]\ar[r]&1\\
1\ar[r]&\Gamma_{d,e}\ar[r]\ar[d]&N_{e,d}\ar[d]\ar[r]&C_n\ar[d]\ar[r]&1\\
1\ar[r]&\Gamma_{\ell,\ell}\ar[r]\ar[d]&N_{\ell,\ell}\ar[d]\ar[r]&C_\ell\ar[d]\ar[r]&1\\
&1&1&1&\\
}
$$ 
with exact rows and columns, 
where $N_{a,b}=C_a*C_b$, and the arrows $N_{d',e'}\to C_m$ and $N_{\ell,\ell}\to C_\ell$
are as in Section \ref{decoprime} and Lemma \ref{schottkygen-q}, respectively.
Thus, $\Gamma_{d',e'}$ and $\Gamma_{\ell,\ell}$ are free of ranks $(d'-1)(e'-1)$ and $\ell-1$, respectively. 
From the diagram, it follows that $\Gamma_{d,e}$ is generated by $\Gamma_{e',d'}$ and
the $C_m$-orbits of $\Gamma_{\ell,\ell}$, where $C_m$ acts by conjugation with the powers
of some primitive $m$-th root of unity contained in  $N_{d,e}$. As the right and the middle columns
are split, also the left column splits. Therefore, $\Gamma_{d,e}$ is free and is of rank
$$
g=(d'-1)(e'-1)+(\ell-1)\cdot m,
$$ 
and we can construct in an obvious way explicit generators for $\Gamma_{d,e}$ from the  generating systems
of $\Gamma_{d',e'}$ and $\Gamma_{\ell,\ell}$ given earlier. Again one checks that these generators
yield a Schottky group if and only if 
$$
\absolute{1-\lambda}<\alpha_p(d,e).
$$
\end{proof}

\begin{rem}
We are convinced that one can refine the method in \cite{vanSteen1982}  in order to relate
to arbitrary $m$-covers
the precise Schottky group, as constructed here.
\end{rem}

\section*{Acknowledgements}

The author thanks his advisor  Frank Herrlich for valuable discussions and support during the
time of writing the dissertation \cite{BradDiss2002}. Especially indebted is the author also to
Stefan K\"uhnlein. It was Fumiharu Kato's encouragement which gave the impetus of writing down this article. Gunther Cornelissen's remarks on a first draft
 and the comments of the unknown referee 
helped improve the exposition of this article.
Finally, Prof.\ Dr.\ Niklaus Kohler's generous wish to see more mathematics really deserves mentioning.

\medskip\noindent
{\sc Universit\"at Karlsruhe, Institut f\"ur Industrielle Bauproduktion, Englerstr.~7, D-76128 Karlsruhe, Germany}\\
E-mail: {\tt bradley@ifib.uni-karlsruhe.de}

\end{document}